\documentclass[11pt]{amsart}

\usepackage[T1]{fontenc}
\usepackage{lmodern}
\usepackage{amsmath,amssymb,amsthm,mathtools}
\usepackage{microtype}
\usepackage{enumitem}
\usepackage[colorlinks=true,linkcolor=blue,citecolor=blue,urlcolor=blue]{hyperref}
\hypersetup{
  pdftitle={Unique Factorization of Power Series Rings over Q-UFDs of Dimension at Most 3},
  pdfauthor={Elad Paran},
  pdfsubject={Unique factorization in formal power series rings over Noetherian UFDs}
}

\setlist[enumerate]{leftmargin=2.5em,itemsep=2pt,topsep=4pt}
\allowdisplaybreaks[1]

\newtheorem{theorem}{Theorem}[section]
\newtheorem{proposition}[theorem]{Proposition}
\newtheorem{lemma}[theorem]{Lemma}
\newtheorem{corollary}[theorem]{Corollary}
\theoremstyle{definition}
\newtheorem{example}[theorem]{Example}
\theoremstyle{remark}
\newtheorem{remark}[theorem]{Remark}
\DeclareMathOperator{\Spec}{Spec}
\DeclareMathOperator{\Pic}{Pic}
\DeclareMathOperator{\Cl}{Cl}
\DeclareMathOperator{\depth}{depth}
\DeclareMathOperator{\Ass}{Ass}
\DeclareMathOperator{\Ann}{Ann}
\DeclareMathOperator{\htop}{ht}
\DeclareMathOperator{\Proj}{Proj}
\newcommand{\m}{\mathfrak m}
\newcommand{\p}{\mathfrak p}
\newcommand{\q}{\mathfrak q}
\newcommand{\cI}{\mathfrak c}

\title[Unique Factorization of Power Series Rings]{Unique Factorization of Power Series Rings over $\mathbb Q$-UFDs of Dimension at Most~3}
\author{Elad Paran}
\address{Department of Mathematics and Computer Science, The Open University of Israel, Ra'anana 4353701, Israel}
\email{paran@openu.ac.il}
\date{}


\begin{document}

\begin{abstract}
Let $R$ be a Noetherian unique factorization domain containing $\mathbb Q$, with $\dim R\leq3$. We prove that $R[[t]]$ is a unique factorization domain if and only if two conditions hold: every two-generated parameter ideal $q$ in every height-two localization satisfies $\overline{q^2}\subseteq q$, and every height-three localization has depth three. No excellence or completion hypothesis is imposed. The main local result is that $A[[t]]$ is not a UFD whenever $A$ is a local $\mathbb Q$-UFD of dimension at least three and depth two. We also prove that the UFD property of $R[[t]]$ passes to every coefficient localization in dimension at most three, give an all-dimensional criterion when the completions of the maximal localizations satisfy Serre's condition $(S_2)$, and obtain an arbitrary-characteristic criterion for excellent UFDs of dimension at most three.
\end{abstract}

\maketitle

\section{Introduction}

Let $A$ be a UFD.\footnote{The older literature often uses the term \emph{factorial ring} for a UFD. We use ``UFD'' throughout.} Is the power series ring $A[[t]]$ a UFD? This problem goes back to Krull, who suspected the general answer is negative \cite{Krull1938}, which was confirmed in 1961 by Samuel \cite[p.~1 and \S4]{SamuelUFD}. Samuel nevertheless conjectured that $A[[t]]$ is a UFD whenever $A$ is a complete local UFD \cite[p.~171]{Samuel1961}. Salmon disproved this conjecture, and Danilov later constructed a systematic family of counterexamples and proposed the property of being geometrically a UFD as a possible replacement \cite[pp.~127--131]{DanilovSamuel}. The UFD property of $A$, even together with completeness, therefore does not determine whether $A[[t]]$ is a UFD.

Danilov reformulated the problem in terms of divisor class groups. A normal Noetherian domain $A$ is said to have a \emph{discrete divisor class group} (DGC) if the canonical map
\[
\Cl(A)\longrightarrow\Cl(A[[t]])
\]
is an isomorphism. For a UFD, the DGC property is equivalent to $A[[t]]$ being a UFD. Danilov proved faithfully flat descent and a local-to-global theorem for DGC, reduced the problem to localizations of depth two, and obtained a characteristic-zero criterion under auxiliary finiteness hypotheses \cite{Danilov}. Lipman removed those hypotheses from the Danilov--Samuel reduction: a normal Noetherian domain has DGC once its depth-two localizations have DGC \cite[Introduction, pp.~203--204]{Lipman1979}.

The local problem has two parts. In dimension two, Lipman related DGC to pseudo-rationality, and Lipman--Teissier related pseudo-rationality to integral-closure identities for ideals \cite[footnote~8, p.~539]{Lipman1975}, \cite[pp.~113--114]{LipmanTeissier}. In dimension at least three, depth two is itself an obstruction: Danilov proved this when the ring is excellent or its completion is normal. His results yield an exact criterion for excellent Noetherian $\mathbb Q$-UFDs in terms of Serre's condition $(S_3)$ and the two-dimensional localizations. The missing case was that of nonexcellent rings, whose completions may be nonnormal, nonreduced, or nonequidimensional.

In this paper, we remove the completion and excellence hypotheses in dimension at most three. First, if $A$ is a local $\mathbb Q$-UFD of dimension at least three and depth two, then $A[[t]]$ is not a UFD, without any assumption on the completion of $A$. Second, at height two we replace pseudo-rationality by a test involving only two-generated parameter ideals, again without assuming that the completion is reduced. Together with the Danilov--Samuel reduction, these results give a criterion expressed entirely in terms of the arithmetic and ideal structure of the coefficient ring. We write $\overline I$ for the integral closure of an ideal $I$.

\begin{theorem}[Main criterion]\label{thm:main}
Let $R$ be a Noetherian UFD containing $\mathbb Q$, and suppose $\dim R\leq3$. Then $R[[t]]$ is a UFD if and only if
\begin{enumerate}[label=\textup{(\roman*)},leftmargin=2.7em]
\item for every prime $\p\subset R$ of height two and all $u,v\in R_\p$ with $\sqrt{(u,v)R_\p}=\p R_\p$,
\begin{equation}\label{eq:surface-test}
\overline{((u,v)R_\p)^2}\ \subseteq\ (u,v)R_\p;
\end{equation}
\item $\depth R_\p=3$ for every prime $\p\subset R$ of height three.
\end{enumerate}
\end{theorem}

For a three-dimensional local UFD, condition (ii) says precisely that $R$ is Cohen--Macaulay. A two-dimensional Noetherian UFD is automatically Cohen--Macaulay, so condition (i) alone applies in that dimension. The qualification ``for all $u,v$'' in (i) is important: the ideal $q=(u,v)R_\p$ is a parameter ideal of the two-dimensional local ring $R_\p$, and the containment is tested in that ring.

Danilov asked whether the DGC property is preserved under localization and proved this for excellent rings containing $\mathbb Q$ \cite[Proposition~4, pp.~371--372]{Danilov}. This was noted as an open problem by Lipman \cite[footnote~2, p.~203]{Lipman1979}. For UFDs, the essential local--global implication concerns only primes $\p$ for which $\depth R_\p=2$. In dimension at most three, the main criterion gives the stronger conclusion at every prime.

\begin{corollary}[Localization of UFD power series rings]\label{cor:localization-intro}
Let $R$ be a Noetherian UFD containing $\mathbb Q$, with $\dim R\leq3$. If $R[[t]]$ is a UFD, then $R_\p[[t]]$ is a UFD for every prime ideal $\p\subset R$.
\end{corollary}
\begin{proof}
By Theorem~\ref{thm:main}, $R$ satisfies the height-two parameter-ideal condition and the height-three depth condition. Both conditions pass to $R_\p$, so another application of the theorem proves the assertion.
\end{proof}

This proves the Danilov--Lipman localization implication for $\mathbb Q$-UFDs of dimension at most three. The corresponding question for arbitrary normal rings with DGC remains open, as does Danilov's conjecture that adjoining a formal variable preserves the property of being geometrically a UFD \cite[p.~131]{DanilovSamuel}.

Lipman's Danilov--Samuel theorem gives the converse implication in every dimension: if $R_\p[[t]]$ is a UFD whenever $\depth R_\p=2$, then $R[[t]]$ is a UFD. Thus an all-dimensional local--global theorem requires coefficient localization only at depth-two primes. Proposition~\ref{prop:conditional-depth-two} states the resulting criterion.

The main local result assumes neither excellence nor reducedness or equidimensionality of the completion.

\begin{theorem}[Depth-two obstruction]\label{thm:local-intro}
If $(A,\m)$ is a Noetherian local UFD containing $\mathbb Q$, $\dim A\geq3$, and $\depth A=2$, then $A[[t]]$ is not a UFD. If $R$ is a Noetherian UFD containing $\mathbb Q$ and a maximal ideal $\m$ satisfies $\htop\m\geq3$ and $\depth R_\m=2$, then $R[[t]]$ is not a UFD.
\end{theorem}

The second statement requires a separate argument because the coefficient localization of $R[[t]]$ is generally smaller than $R_\m[[t]]$. We construct a divisor class that descends from $\widehat{R_\m}[[t]]$ to $R[[t]]$.

Theorem~\ref{thm:main} extends the excellent-ring classification to all Noetherian $\mathbb Q$-UFDs of dimension at most three. Section~\ref{sec:higher-dimensions} also gives a criterion in any dimension when maximal completions satisfy $(S_2)$, and identifies the case left open in dimension four.

Throughout the paper, all rings are commutative and Noetherian. The ring $A[[t]]$ consists of formal power series with coefficients in $A$. If $A$ is a domain, $\Cl(A)$ denotes its Weil divisor class group, represented by rank-one reflexive modules. We write $H^i_I(M)$ for local cohomology with support in $V(I)$.

Section~2 proves the surface criterion and its globalization. Sections~3 and~4 develop the deformation and completion arguments used in the depth-two obstruction. Section~5 proves the main theorem. Section~6 treats higher dimensions, and Section~7 gives an arbitrary-characteristic criterion for excellent UFDs of dimension at most three using Witt-vector cohomology. The latter criterion refers only to localizations of the coefficient ring.

\section{The surface test and its globalization}

We begin by extending the Lipman--Teissier surface criterion to rings whose completions may be nonreduced.

\begin{proposition}[Surface parameter ideals]\label{prop:surface}
Let $(A,\m)$ be a two-dimensional normal Noetherian local domain. The following are equivalent:
\begin{enumerate}[label=\textup{(\alph*)},leftmargin=2.7em]
\item $A$ is pseudo-rational;
\item $\overline{q^2}\subseteq q$ for every two-generated $\m$-primary ideal $q\subset A$.
\end{enumerate}
If $A$ is a UFD, these conditions are also equivalent to $A[[t]]$ being a UFD. The assertion holds in every characteristic and does not assume excellence.
\end{proposition}
\begin{proof}
If $A$ is pseudo-rational, Lipman--Teissier's Corollary~5.4 gives $\overline{I^2}=I\overline I\subseteq I$ for every ideal $I$, proving (b). Conversely, suppose (b). First, $\widehat A$ is reduced. Indeed, if $0\ne z\in\widehat A$ and $z^e=0$, fix parameters $u,v\in A$. For some $N$, the ideal $q=(u^N,v^N)A$ satisfies $z\notin q\widehat A$: the ideals $(u^N,v^N)\widehat A$ shrink inside powers of $\m\widehat A$ and their intersection is zero. Choose $M$ so that $\m^M\widehat A\subseteq q^{2e}\widehat A$, and choose $a\in A$ with $a-z\in\m^M\widehat A$. The binomial expansion gives $a^e\in q^{2e}\widehat A\cap A=q^{2e}$, since completion is faithfully flat. Hence $a\in\overline{q^2}$, whereas $a\notin q$ because $z\notin q\widehat A$. This contradicts (b).

Every parameter pair in $A$ is a regular sequence. The Huneke--Itoh intersection theorem \cite[Theorem~1]{Itoh} therefore gives $q\cap\overline{q^2}=q\overline q$. By (b), $\overline{q^2}=q\overline q$ for all parameter ideals $q$. Now $\widehat A$ is reduced, so the converse argument of Lipman--Teissier's Proposition~5.5 applies: their Corollary~5.3(i) computes the first cohomology on normalized blowups from $\overline{q_n^2}/(q_n\overline{q_n})$, for parameter ideals $q_n=(u^n,v^n)$, and those blowups dominate the normal proper birational models. It follows that $A$ is pseudo-rational \cite[pp.~113--114]{LipmanTeissier}. The last equivalence is Lipman's two-dimensional class-group result \cite[footnote~8, p.~539]{Lipman1975}.
\end{proof}
The assumption that $R[[t]]$ is a UFD itself forces the surface condition; this does not follow formally by localizing $R[[t]]$.

\begin{proposition}[Global surface obstruction]\label{prop:height-two}
Let $R$ be a Noetherian UFD. If $R[[t]]$ is a UFD, then \eqref{eq:surface-test} holds at every height-two prime $\p$. Conversely, if $R$ satisfies $(S_3)$ and \eqref{eq:surface-test} holds at every height-two prime, then $R[[t]]$ is a UFD. Neither assertion needs $\mathbb Q$ or excellence.
\end{proposition}
\begin{proof}
For necessity, suppose $A=R_\p$, $q=(u,v)A$, and $z\in\overline{q^2}\setminus q$. Clearing the finitely many denominators of $u,v,z$ and of an equation of integral dependence, and multiplying $z$ by an element outside $\p$, gives $u,v,z\in R$ with $z\in\overline{((u,v)R)^2}$ and $z\notin(u,v)A$. Let $s$ be a new variable and put
\[
C=R[us,vs,zs^2]\subset R[s],\qquad X=\Proj C.
\]
The element $zs^2$ is integral over the Rees algebra $R[us,vs]$. Thus $X$ is projective over $R$, and the two standard opens $D_+(us)$ and $D_+(vs)$ cover it. The regular function $h=z/(uv)$ is defined on their intersection. Put $B=R[[t]]$ and $X_n=X\times_R\Spec(R[t]/t^n)$. Gluing the trivial line bundle on the two charts with transition $1+th$ defines compatible line bundles $\mathcal L_n$ on $X_n$. The inverse exists modulo $t^n$. Grothendieck existence for the projective $B$-scheme $X\times_R\Spec B$, with $B$ complete for $(t)$, algebraizes them to a line bundle $\mathcal L$ \cite[Theorem~30.27.1, Tag~088E]{Stacks}. To see that the algebraization is invertible, apply existence also to the inverse bundles and use full faithfulness for the evaluation maps.

Over $W=D(u)\cup D(v)\subset\Spec B$, this projective scheme equals $W$: after inverting $u$ or $v$, its graded algebra is a polynomial ring in $s$. If $B$ were a UFD, then $\Pic(W)=0$ \cite[Lemma~31.29.4, Tag~0BDA]{Stacks}. Thus the restriction of $\mathcal L_2$ to the punctured spectrum of $A[t]/t^2$ would be trivial. On that spectrum its transition is $1+t z/(uv)$. The split square-zero unit sequence identifies the kernel of reduction $\Pic(U\times_A\Spec(A[t]/t^2))\to\Pic(U)$ with $tH^1(U,\mathcal O_U)$, where $U=D(u)\cup D(v)=\Spec A\setminus\{\p A\}$. Here
\[
H^1(U,\mathcal O_U)=A_{uv}/(A_u+A_v).
\]
Moreover, the map
\[
A/(u,v)A\longrightarrow A_{uv}/(A_u+A_v),
\qquad \bar a\longmapsto[a/(uv)],
\]
is injective because $u,v$ is an $A$-regular sequence. The class of $z/(uv)$ is therefore nonzero, a contradiction.

For sufficiency, $(S_3)$ implies that every localization of $R$ of depth two has height two. Proposition~\ref{prop:surface} makes $R_\p[[t]]$ a UFD at every such prime. The Danilov--Samuel theorem in Lipman's unrestricted form \cite[Introduction, pp.~203--204]{Lipman1979} yields $\Cl(R)\simeq\Cl(R[[t]])$. Both groups vanish because $R$ is a UFD.
\end{proof}

\begin{corollary}[The local--global theorem for $(S_3)$ UFDs]
\label{cor:s3-local-global}
Let $R$ be a Noetherian UFD satisfying $(S_3)$. If $R[[t]]$ is a UFD, then $R_\p[[t]]$ is a UFD for every prime $\p\subset R$.
\end{corollary}
\begin{proof}
Fix $\p$ and put $A=R_\p$. If $\q A$ is a prime of $A$ such that $\depth A_{\q A}=2$, then $A_{\q A}=R_\q$ and $(S_3)$ forces $\htop\q=2$. Proposition~\ref{prop:height-two}, followed by Proposition~\ref{prop:surface}, shows that $R_\q[[t]]$ is a UFD. Lipman's Danilov--Samuel theorem now gives that $A[[t]]$ is a UFD.
\end{proof}

\begin{example}[Why the surface condition is needed]\label{ex:danilov}
For Danilov's complete local UFD $A_0=\mathbb Q[[x,y,z]]/(2x^3+4y^3-3z^3)$, the ring $A_0[[t]]$ is not a UFD \cite[\S1, p.~128]{DanilovSamuel}. Indeed $q=(x,y)A_0$ is $\m$-primary, $z^2\notin q$ since $A_0/q=\mathbb Q[[z]]/(z^3)$, and $(z^2)^3=z^6=((2x^3+4y^3)/3)^2\in q^6=(q^2)^3$. Thus $z^2\in\overline{q^2}\setminus q$.
\end{example}

\section{Three ways to produce a divisor class}

We develop the obstructions for a complete local ring $T$. They will also be used when $T$ is the completion of a localization of a global ring. The first construction works before completion.

\begin{proposition}[Finite local cohomology]\label{prop:finite-H2}
Let $R$ be a Noetherian UFD containing $\mathbb Q$ and let $J\subset R$ have height at least two. If $H^2_J(R)$ is a nonzero finite $R$-module, then $R[[t]]$ is not a UFD.
\end{proposition}
\begin{proof}
Set $B=R[[t]]$, $U=\Spec R\setminus V(J)$, and $X=\Spec B\setminus V(JB+(t))$. Since $R$ is normal and $\htop J\geq2$, $H^1_J(R)=0$, $\Gamma(U,\mathcal O_U)=R$, and $H^1(U,\mathcal O_U)=H^2_J(R)$. Choose $0\ne\alpha\in H^1(U,\mathcal O_U)$ and a Cech representative $\alpha_{ij}$ on a finite principal-open cover. On $X_n=X\times_B\Spec(B/t^nB)$ the expressions
\[
\exp(t\alpha_{ij})\pmod{t^n}=\sum_{a=0}^{n-1}(t\alpha_{ij})^a/a!
\]
are a multiplicative cocycle. They define compatible line bundles $L_n$, with $L_1$ trivial. The split square-zero unit sequence gives $[L_2]=t\alpha\ne0$ in the kernel of $\Pic(X_2)\to\Pic(U)$.

The groups $H^1_{JB+(t)}(B/tB)=H^1_J(R)=0$ and $H^2_{JB+(t)}(B/tB)=H^2_J(R)$ are finite. By the formal algebraization theorem \cite[Lemma~52.16.11, Tag~0DXW]{Stacks}, there is a coherent sheaf $\mathcal F$ on $X$ with completions $L_n$. It is locally free of rank one on an open $V\subseteq X$ containing the entire special fiber $U$: at a point of $U$, completion along $t$ is free, and completion is faithfully flat on the local ring. If $B$ were a UFD then $\Pic(V)=0$ \cite[Lemma~31.29.4, Tag~0BDA]{Stacks}, making $L_2$ trivial. This contradicts $[L_2]=t\alpha\ne0$.
\end{proof}

The following descent lemma is stated over a general Noetherian base so that it applies at a maximal ideal of a nonlocal ring.

\begin{lemma}[Descent with finite self-Ext]\label{lem:descent}
Let $R$ be a Noetherian $\mathbb Q$-algebra, $\m$ a maximal ideal, $A=R_\m$, and $T=\widehat A$. Let $J\subset R$ be an $\m$-primary ideal. Suppose $N$ is a finite $T[[t]]$-module, multiplication by $t$ is injective on $N$, and $N/tN\simeq JT$. Then there is a finite $R[[t]]$-module $M$ such that $M/tM\simeq J$ and
\begin{equation}\label{eq:descent}
M\otimes_{R[[t]]}T[[t]]\simeq N.
\end{equation}
The same assertion holds for a finite $R$-module $J$ if $\operatorname{Ext}^i_R(J,J)$ has finite length supported at $\m$ for $i=1,2$.
\end{lemma}
\begin{proof}
The maps $R\to A\to T$ identify $\operatorname{Ext}^i_R(J,J)$ with $\operatorname{Ext}^i_T(JT,JT)$ for $i=1,2$. Indeed, outside $\m$ the ideal $J$ is $R$, so its positive self-Ext groups vanish; being finite $R$-modules supported at $\m$, they have finite length, and completion does not change them.

Write $R_n=R[t]/(t^n)$, $T_n=T[t]/(t^n)$, and $N_n=N/t^nN$. Because $t$ is regular on $N$ and $\mathbb Q\subseteq T$, $N_n$ is flat over $\mathbb Q[t]/(t^n)$. Starting with $M_1=J$ and the given isomorphism $M_1\otimes_RT\simeq N_1$, construct compatible $R_n$-modules $M_n$, flat over $\mathbb Q[t]/(t^n)$, together with framed isomorphisms $M_n\otimes_RT\simeq N_n$.

For the step $n\to n+1$, the obstruction and the torsor of framed lifts lie in $\operatorname{Ext}^2_{R_n}(M_n,J)$ and $\operatorname{Ext}^1_{R_n}(M_n,J)$, respectively \cite[Tags~08VR and 08LM]{Stacks}. Flatness over $\mathbb Q[t]/(t^n)$ gives $M_n\otimes^{\mathbf L}_{R_n}R\simeq J$, so derived adjunction identifies these Ext groups with $\operatorname{Ext}^i_R(J,J)$. Their maps to the corresponding groups over $T$ are isomorphisms. The obstruction vanishes over $T$ because $N_{n+1}$ exists; it therefore vanishes over $R$. The isomorphism on Ext$^1$ permits the chosen framed lift to match $N_{n+1}$. Affine formal effectivity \cite[Lemma~30.23.1, Tag~087W]{Stacks} yields $M=\varprojlim M_n$, a finite module over $R[[t]]$. The compatible identifications modulo every $t^n$ yield \eqref{eq:descent}.
\end{proof}

\begin{remark}\label{rem:basechange}
The map $R[[t]]\to T[[t]]$ is flat: it is the composite of the $t$-adic completion of $(R\setminus\m)^{-1}R[[t]]$ and the $\m$-adic completion of $R_\m[[t]]$. Flat base change therefore commutes with the double dual of a finite module. If $R[[t]]$ is a UFD, every rank-one reflexive $R[[t]]$-module is free. We use these facts after applying Lemma~\ref{lem:descent}; they do not assert that $R_\m[[t]]$ itself is a localization of $R[[t]]$.
\end{remark}

\begin{proposition}[Conductor obstruction]\label{prop:conductor}
Let $(T,\m_T)$ be a complete Noetherian local ring containing $\mathbb Q$. Suppose $T\hookrightarrow S$ is a finite extension of rings and $C=S/T$ satisfies
\[
\depth_{\m_T}S\geq2,\qquad C\ne0,\qquad \dim_T C=\depth_{\m_T}C=1.
\]
If $R$ is a Noetherian UFD containing $\mathbb Q$ with a maximal ideal $\m$ such that $\widehat{R_\m}\simeq T$, then $R[[t]]$ is not a UFD. This includes $R=A$ for every local UFD $A$ with completion $T$.
\end{proposition}
\begin{proof}
The local cohomology sequence of $0\to T\to S\to C\to0$ gives $\depth T=2$: it gives $H^0_{\m_T}(T)=H^1_{\m_T}(T)=0$ and injects the nonzero $H^1_{\m_T}(C)$ into $H^2_{\m_T}(T)$.

Let $\cI=(T:S)=\Ann_T C$ and put $H=T/\cI$, $K=S/\cI$. The conductor $\cI$ is an ideal of $S$. Because $C$ is faithful as an $H$-module, $\dim H=1$. Its depth is one: an element in the socle of $H$ would act nontrivially on $C$ and produce a nonzero element killed by $\m_T$, contrary to $\depth C=1$. From $0\to H\to K\to C\to0$ it follows that $K$ also has depth and dimension one. Choose $r\in\m_T$ regular on $H$, $K$, and $C$, and choose $\xi\in K$ whose image $\bar\xi\in C/rC$ is nonzero. Write $E=T[[t]]$ and $S'=S[[t]]$. Define a finite $E$-submodule $N\subseteq S'$ by
\begin{equation}\label{eq:conductor-N}
N=\{s\in S':s\bmod\cI S'\in(r+t\xi)H[[t]]\subset K[[t]]\}.
\end{equation}
It is $t$-torsionfree. Moreover $N\cap tS'=tN$: if $ts\in N$, the constant coefficient of $(r+t\xi)h=ts\bmod\cI$ is $rh_0=0$ in $K$, so $h_0=0$ and one can cancel $t$. Consequently
\[
N/tN\simeq (r,\cI)T=:I.
\]
Since $H/rH$ has finite length, $I$ is $\m_T$-primary. For every nonmaximal prime $\q$ of $T$, $I_\q=T_\q$; regularity of $t$ and Nakayama's lemma make $N$ locally free of rank one at the points of the special fiber above $\q$. In particular $N=N^{**}$ at these points.

Put $U=\Spec T\setminus\{\m_T\}$. The opens $D(r)$ and $D(c_i)$, for generators $c_i$ of $\cI$, cover $U$. On the latter opens $N/t^2N$ has generator $1$; on $D(r)$ it has generator $r+t\widetilde\xi$, for any lift $\widetilde\xi\in S$ of $\xi$. Because $c_iS\subset T$, on intersections $D(r c_i)$ the first-order transition is $1+t\widetilde\xi/r$. In the local cohomology sequence for $0\to T\to S\to C\to0$, $\depth_{\m_T}S\geq2$ and $\Gamma(U,\widetilde C)=C_r$ give an injection
\begin{equation}\label{eq:conductor-class}
C_r/C\ \hookrightarrow\ H^1(U,\mathcal O_U)=H^2_{\m_T}(T).
\end{equation}
The class of the first-order line bundle corresponds to $\bar\xi/r$ under this injection. It is nonzero: if $\bar\xi/r\in C$, then $\bar\xi\in rC$. Hence $N/t^2N$ is a nontrivial line bundle on $U\times_T\Spec(T[t]/t^2)$, trivialized at $t=0$.

Let $J=I\cap R$. Since $I$ contains a power of $\m_T$, the isomorphisms $R/\m^a\simeq T/\m_T^a$ give $JT=I$ and show that $J$ is $\m$-primary. Lemma~\ref{lem:descent} yields $M$ over $R[[t]]$ whose completion is $N$ and whose special fiber is $J$. The module $M$ has generic rank one: at the height-one prime $(t)$, its special fiber is the fraction field of $R$, and $t$ is regular on $M$. If $R[[t]]$ were a UFD, $M^{**}$ would be free of rank one. Remark~\ref{rem:basechange} would make $N^{**}$ free, which would trivialize the line bundle $N/t^2N$ on $U\times_T\Spec(T[t]/t^2)$ because $N=N^{**}$ there. This contradicts \eqref{eq:conductor-class}.
\end{proof}

The nilpotent case requires an exponential with only finitely many terms; no formal algebraization is needed for the line bundle itself.

\begin{proposition}[Nilpotent obstruction]\label{prop:nilpotent}
Let $(T,\m_T)$ be a complete Noetherian local $\mathbb Q$-algebra of depth at least two. If $P\subseteq T$ is a nilpotent ideal such that
\begin{equation}\label{eq:nil-map}
\operatorname{im}\bigl(H^2_{\m_T}(P)\longrightarrow H^2_{\m_T}(T)\bigr)\ne0,
\end{equation}
then $R[[t]]$ is not a UFD for any Noetherian UFD $R\supseteq\mathbb Q$ with a maximal ideal $\m$ and $\widehat{R_\m}\simeq T$.
\end{proposition}
\begin{proof}
Let $U=\Spec T\setminus\{\m_T\}$, and choose $\alpha\in H^1(U,\widetilde P)$ whose image $\beta\in H^1(U,\mathcal O_U)$ is nonzero. Put $E=T[[t]]$ and $W=\Spec E\setminus V(\m_TE)$. Represent $\alpha$ by a Cech cocycle $\alpha_{ij}$ for a finite principal-open cover of $U$. On the corresponding cover of $W$, $g_{ij}=\exp(t\alpha_{ij})$ is a unit and the $g_{ij}$ form a multiplicative cocycle: $P$ is nilpotent, so every exponential is a finite sum. Let $\mathcal L$ be the resulting line bundle. Its reduction to the first infinitesimal thickening of $U$ has class $t\beta\ne0$, hence $\mathcal L$ is nontrivial.

The module $N=\Gamma(W,\mathcal L)$ is finite over $E$. To apply the coherence criterion \cite[Lemma~51.8.3, Tag~0AWA]{Stacks}, note that $E$ is a complete Noetherian local ring, thus Nagata and universally catenary. Its associated primes are $\p E$ for $\p\in\Ass T$; every such $\p$ satisfies $\dim T/\p\geq\depth T\geq2$ \cite[Lemma~10.72.9, Tag~0BK4]{Stacks}. The two points of $V(\m_TE)$ have codimension at least two in every such associated component. Thus the pushforward of $\mathcal L$ is coherent. Applying the same argument to $\mathcal L^{-1}$ and using $\Gamma(W,\mathcal O_W)=E$ shows $N=N^{**}$; it is rank one on every minimal component and $t$ acts injectively on it.

Reduction of sections gives an injection $N/tN\hookrightarrow\Gamma(U,\mathcal O_U)=T$: a section that reduces to zero can be divided by $t$ locally, and its quotients glue. Call its image $I$. Since $\mathcal L/t\mathcal L=\mathcal O_U$, one has $I_\q=T_\q$ for every nonmaximal prime $\q$ of $T$. Thus $I=T$ or $I$ is $\m_T$-primary. The case $I=T$ would make $N$ free by Nakayama's lemma, contrary to the nontriviality of $\mathcal L$. Therefore $I$ is $\m_T$-primary.

Contract $I$ to an $\m$-primary ideal $J=I\cap R$. The quotients $R/\m^a\simeq T/\m_T^a$ imply $JT=I$. Descend $N$ by Lemma~\ref{lem:descent}. If $R[[t]]$ were a UFD, the double dual of the descended rank-one module would be free; its flat base change $N$ would then be free, contradicting $[\mathcal L]\ne0$.
\end{proof}

\section{Arbitrary depth-two completions}

We first treat completions whose associated components all have dimension at least three. In the next proof, ``finite'' means finitely generated as a module rather than finite length.

\begin{lemma}\label{lem:large-components}
Let $(T,\m_T)$ be a complete Noetherian local $\mathbb Q$-algebra with $\depth T=2$ and $\dim T/\p\geq3$ for every $\p\in\Ass T$. For every Noetherian UFD $R\supseteq\mathbb Q$ with a maximal ideal $\m$ and $\widehat{R_\m}\simeq T$, the ring $R[[t]]$ is not a UFD.
\end{lemma}
\begin{proof}
If $H^2_{\m_T}(T)$ is finite, it is nonzero because $\depth T=2$. Flat base change, faithful flatness of $R_\m\to T$, and localization at $\m$ show that $H^2_\m(R)$ is a nonzero finite $R$-module. Proposition~\ref{prop:finite-H2}, with $J=\m$, proves the assertion.

Suppose $H^2_{\m_T}(T)$ is not finite. The finiteness criterion of SGA~2 \cite[Expos\'e~VIII, Corollaire~2.3]{SGA2} gives a prime $\q\subsetneq\m_T$ with $\dim(T/\q)=1$ and $\depth T_\q<2$. This applies because $T$ is a quotient of a regular complete local ring and $H^0_{\m_T}(T)=H^1_{\m_T}(T)=0$. Our hypothesis excludes $\q\in\Ass T$; hence $\depth T_\q=1$.

Set $U=\Spec T\setminus V(\q)$ and $S=\Gamma(U,\mathcal O_U)$. The missing set is $\{\q,\m_T\}$. For every $\p\in\Ass T$ lying below $\q$, catenarity of $T/\p$ gives
\[
\operatorname{codim}_{V(\p)}(\q)=\dim(T/\p)-1\geq2,
\]
and the closed point also has codimension at least two in every associated component. Since $T$ is Nagata and universally catenary, the coherent-pushforward criterion \cite[Lemma~51.8.3, Tag~0AWA]{Stacks} gives that $S$ is finite over $T$. The natural map $T\to S$ is injective: its kernel would have an associated prime containing $\q$, which would contradict $\dim T/\p\geq3$ for all $\p\in\Ass T$.

Let $C=S/T$. Its support lies in $V(\q)$, and localization at $\q$ gives
\[
C_\q\simeq H^1_{\q T_\q}(T_\q)\ne0.
\]
Thus $\dim C=1$. No nonzero section of $S=\Gamma(U,\mathcal O_U)$ can be killed by a power of $\m_T$; hence $H^0_{\m_T}(S)=0$. Since $H^1_{\m_T}(T)=0$, the sequence $0\to T\to S\to C\to0$ yields $H^0_{\m_T}(C)=0$. Therefore $\depth_{\m_T}C=1$.

Finally, write $W=\Spec T\setminus\{\m_T\}$. The sheaf associated to $S$ is $j_*\mathcal O_U$ for $j:U\hookrightarrow\Spec T$. Thus $\Gamma(W,\widetilde S)=\Gamma(U,\mathcal O_U)=S$. The local cohomology sequence of $S$ then gives $H^1_{\m_T}(S)=0$. We have produced a finite extension $T\hookrightarrow S$ satisfying Proposition~\ref{prop:conductor}, which completes the proof.
\end{proof}

\begin{theorem}[Complete local depth-two obstruction]\label{thm:completion}
Let $(T,\m_T)$ be any complete Noetherian local ring containing $\mathbb Q$ with $\dim T\geq3$ and $\depth T=2$. If $R$ is a Noetherian UFD containing $\mathbb Q$ and $\m$ is a maximal ideal with $\widehat{R_\m}\simeq T$, then $R[[t]]$ is not a UFD. Every such $T$ occurs as the completion of a Noetherian local UFD.
\end{theorem}
\begin{proof}
Heitmann's completion theorem \cite[Theorems~1 and~8]{Heitmann} gives the last assertion: the depth is two and every nonzero integer is a unit. Fix any $R$ and $\m$ as in the statement. The associated-prime depth inequality \cite[Lemma~10.72.9, Tag~0BK4]{Stacks} says $\dim(T/\p)\geq2$ for all $\p\in\Ass T$. Lemma~\ref{lem:large-components} handles the case where all these dimensions are at least three. Suppose, therefore, that some $\p\in\Ass T$ has $\dim(T/\p)=2$.

Write $P=\sqrt{(0)}\subset T$ and $D=T/P$. The reduced complete local ring $D$ has no minimal component of dimension below two. Its normalization $\overline D$ is a finite product of complete normal local domains, each of dimension at least two, and $\depth_{\m_T}\overline D\geq2$. The exact sequence $0\to D\to\overline D\to\overline D/D\to0$ therefore gives
\begin{equation}\label{eq:H1-reduced}
H^1_{\m_T}(D)\simeq H^0_{\m_T}(\overline D/D),
\end{equation}
a finite-length module.

\emph{Case 1: $P$ has a two-dimensional associated component.} Suppose $\p\in\Ass_T P$ and $\dim T/\p=2$. Choose a prime $\q\supsetneq\p$ with $\dim T/\q=1$ and $\htop(\q/\p)=1$ in the complete two-dimensional local domain $T/\p$. The associated prime $\p T_\q$ of $P_\q$ gives $\depth P_\q\leq1$. Also $H^0_{\m_T}(P)=H^1_{\m_T}(P)=0$: the first follows from $P\subset T$ and $\depth T=2$, and the second from $0\to P\to T\to D\to0$ and $H^0_{\m_T}(D)=0$. The SGA~2 finiteness criterion \cite[Expos\'e~VIII, Corollaire~2.3]{SGA2} now shows $H^2_{\m_T}(P)$ is not finite. By \eqref{eq:H1-reduced}, the kernel of $H^2_{\m_T}(P)\to H^2_{\m_T}(T)$ is an image of a finite-length module. The image is therefore nonzero. Proposition~\ref{prop:nilpotent} gives the result.

\emph{Case 2: every two-dimensional associated component is reduced.} Suppose no $\p\in\Ass_T P$ has $\dim T/\p=2$. Because $\Ass T\subseteq\Ass P\cup\Ass D$ and the associated primes of a reduced Noetherian ring are minimal, our chosen $\p\in\Ass T$ with $\dim T/\p=2$ is minimal in $T$. Choose an irredundant primary decomposition
\[
0=Q_\p\cap J,\qquad Q_\p\text{ is $\p$-primary},\qquad
J=\bigcap_{i\ne\p}Q_i.
\]
There is at least one other component since $\dim T\geq3$. Set $B_1=T/Q_\p$ and $B_2=T/J$. The map $T\hookrightarrow B_1\times B_2$ is injective. Each associated prime of each $B_i$ has quotient dimension at least two. For $i=1,2$, set
\[
S_i=\Gamma(\Spec B_i\setminus\{\m_{B_i}\},\mathcal O).
\]
Each $B_i$ is complete and has codimension at least two at its closed point in each associated component. Again \cite[Lemma~51.8.3, Tag~0AWA]{Stacks} makes $S_i$ finite over $B_i$. The maps $B_i\hookrightarrow S_i$ are injective because the closed point is not associated. By their definition as punctured-spectrum sections, $H^0_{\m_{B_i}}(S_i)=H^1_{\m_{B_i}}(S_i)=0$. Thus $S=S_1\times S_2$ is a finite $T$-algebra containing $T$ with $\depth_{\m_T}S\geq2$.

The quotient $(B_1\times B_2)/T$ is $T/(Q_\p+J)$ and has dimension at most one. Indeed, every other associated prime of $T$ is incomparable with $\p$: a prime strictly above $\p$ would have quotient dimension at most one, contrary to $\depth T=2$. The quotients $S_i/B_i$ are supported at their closed points. Hence $C=S/T$ has dimension at most one. It is nonzero, since a local ring cannot be the nontrivial product $B_1\times B_2$. Finally $H^1_{\m_T}(T)=0$ and $H^0_{\m_T}(S)=0$ imply $H^0_{\m_T}(C)=0$. It follows that $C$ is Cohen--Macaulay of dimension one. Proposition~\ref{prop:conductor} completes the proof.
\end{proof}

\section{The dimension-three criterion}

\begin{proof}[Proof of Theorem~\ref{thm:local-intro}]
For the local assertion, take $T=\widehat A$. Completion preserves dimension and depth, so Theorem~\ref{thm:completion} applied with $R=A$ proves it. For the global assertion, take $T=\widehat{R_\m}$; Theorem~\ref{thm:completion} applies directly to $R$. In particular, this step uses global descent in Lemma~\ref{lem:descent} or Proposition~\ref{prop:finite-H2}, and never infers that $R[[t]]$ is a UFD from a coefficient-localized series ring.
\end{proof}

\begin{proof}[Proof of Theorem~\ref{thm:main}]
If $R[[t]]$ is a UFD, Proposition~\ref{prop:height-two} gives (i). A prime $\p$ of height three is maximal because $\dim R\leq3$. Since $R$ is normal, $\depth R_\p\geq2$; if it were two, Theorem~\ref{thm:local-intro} would contradict the assumption that $R[[t]]$ is a UFD. Thus $\depth R_\p=3$, proving (ii).

Conversely, every UFD satisfies $(S_2)$. Condition (ii) therefore implies $(S_3)$ at every prime of $R$, because the ring has dimension at most three. Proposition~\ref{prop:height-two} combines $(S_3)$ with (i) to show $R[[t]]$ is a UFD.
\end{proof}

\begin{example}[A completion with two components]\label{ex:mixed}
Let
\[
T=\mathbb Q[[x,y,z,w]]/(xz,yz).
\]
It is a complete reduced local ring of dimension three and depth two. Its components have dimensions three and two, and its normalization is $\mathbb Q[[x,y,w]]\times\mathbb Q[[z,w]]$. Heitmann's theorem gives a three-dimensional local UFD $A\supseteq\mathbb Q$ with $\widehat A\simeq T$. Such $A$ cannot be excellent, since the completion of an excellent local domain is a domain. Theorem~\ref{thm:completion} states that $A[[t]]$ is not a UFD for \emph{every} local UFD $A$ with this completion. More concretely, the normalization $S$ has $S/T\simeq\mathbb Q[[w]]$, so Proposition~\ref{prop:conductor} already applies.
\end{example}

\subsection{Relation with classical results}

We now describe the precise relation between the preceding results and the classical criteria. Below, $(S_3)$ denotes Serre's condition $(S_3)$: $\depth A_\q\geq\min\{3,\htop\q\}$ for every prime $\q\subseteq A$.
\begin{enumerate}[label=\textup{(\arabic*)},leftmargin=2.65em,itemsep=.35em]
\item \textbf{Excellent $\mathbb Q$-UFDs.} Let $A$ be an excellent Noetherian UFD containing $\mathbb Q$. Then $A[[t]]$ is a UFD if and only if $A$ satisfies $(S_3)$ and, for every height-two prime $\p\subset A$ and every $\p A_\p$-primary ideal $I\subset A_\p$, the equality $\overline{I^2}=I\overline I$ holds in $A_\p$. This is the UFD specialization of Danilov's Theorem~2 \cite[Theorem~2, p.~382]{Danilov}.
\item \textbf{Danilov's local depth obstruction.} Let $(A,\m)$ be a Noetherian local UFD containing $\mathbb Q$, with $\dim A\geq3$ and $\depth A=2$. If $A$ is excellent or its completion $\widehat A$ is normal, then $A[[t]]$ is not a UFD \cite[Proposition~5 and its corollary, pp.~375--376, using Lemma~3, pp.~373--374]{Danilov}.
\item \textbf{The Danilov--Samuel reduction.} Let $A$ be a Noetherian UFD, in any characteristic. If $A$ satisfies $(S_3)$ and $A_\p$ is pseudo-rational for every height-two prime $\p\subset A$, then $A[[t]]$ is a UFD. No excellence assumption is needed \cite[(D--S), Introduction, pp.~203--204]{Lipman1979}, \cite[footnote~8, p.~539]{Lipman1975}.
\item \textbf{The surface criteria of Lipman and Lipman--Teissier.} Let $(A,\m)$ be a two-dimensional Noetherian local UFD. Then $A[[t]]$ is a UFD if and only if $A$ is pseudo-rational \cite[footnote~8, p.~539]{Lipman1975}. If, in addition, $\widehat A$ is reduced, then $A[[t]]$ is a UFD if and only if $\overline{I^2}=I\overline I$ for every $\m$-primary ideal $I$ of $A$. This last equivalence combines Lipman's result with Lipman--Teissier's Corollary~5.4 and Proposition~5.5 \cite[pp.~113--114]{LipmanTeissier}.
\end{enumerate}

We spell out the translations. For item~(1), put $B=A_\p$ at a height-two prime. Danilov's Theorem~2 asks for $(S_3)$ and finiteness of the divisor class group of the strict henselization of $B$. His Proposition~8 identifies the latter condition with preservation of the divisor class group from $B$ to $B[[t]]$ \cite[pp.~380--382]{Danilov}. Since $B$ is a UFD, this is equivalent to $B[[t]]$ being a UFD, which Lipman identifies with pseudo-rationality of $B$ \cite[footnote~8, p.~539]{Lipman1975}. The completion of an excellent normal ring is normal \cite[Lemma~15.53.6, Tag~0C23]{Stacks}, and Lipman--Teissier identify pseudo-rationality with the ideal equalities in item~(1) \cite[Corollary~5.4 and Proposition~5.5, pp.~113--114]{LipmanTeissier}.

For item~(3), Lipman's hypothesis requires $A_\p[[t]]$ to be a UFD whenever $\depth A_\p=2$ \cite[(D--S), Introduction, pp.~203--204]{Lipman1979}. Under $(S_3)$ these primes have height two, and item~(4) identifies the hypothesis with pseudo-rationality of every height-two localization. Without $(S_3)$ one must also treat primes of larger height and depth two. Theorem~\ref{thm:local-intro} rules them out over $\mathbb Q$. Proposition~\ref{prop:surface} further replaces pseudo-rationality, without excellence, by the explicit containment $\overline{q^2}\subseteq q$ for every two-generated parameter ideal $q$.

Thus Proposition~\ref{prop:surface} sharpens item~(4) to the parameter-ideal test \eqref{eq:surface-test}, including rings with nonreduced completions, while Theorem~\ref{thm:local-intro} removes the extra assumptions in item~(2), even for maximal localizations of nonlocal UFDs. Heitmann's theorem permits local $\mathbb Q$-UFDs whose completions have nilpotents or components of different dimensions \cite[Theorems~1 and~8]{Heitmann}.

\section{Dimension four and higher}\label{sec:higher-dimensions}

The surface condition \eqref{eq:surface-test} and the depth obstruction give the following bounds in every dimension.

\begin{corollary}[Bounds in arbitrary dimension]\label{cor:all-dimensions}
Let $R$ be a Noetherian UFD containing $\mathbb Q$, of arbitrary finite dimension. Write $(\mathrm H_2)$ for condition \eqref{eq:surface-test} at every height-two prime. Then:
\begin{enumerate}[label=\textup{(\alph*)},leftmargin=2.7em]
\item If $R[[t]]$ is a UFD, then $(\mathrm H_2)$ holds and $\depth R_\m\geq3$ at every maximal ideal $\m$ of height at least three.
\item If $(\mathrm H_2)$ holds and $R$ satisfies $(S_3)$, then $R[[t]]$ is a UFD.
\item If $R$ is excellent, then $R[[t]]$ is a UFD if and only if $R$ satisfies $(S_3)$ and $(\mathrm H_2)$ holds.
\end{enumerate}
\end{corollary}
\begin{proof}
Part~(a) is Proposition~\ref{prop:height-two} together with Theorem~\ref{thm:local-intro}; part~(b) is the converse in Proposition~\ref{prop:height-two}. For part~(c), Danilov's Theorem~2 and Proposition~8 \cite[pp.~380--382]{Danilov}, combined with Lipman's surface criterion \cite[footnote~8, p.~539]{Lipman1975}, identify Danilov's height-two hypothesis with pseudo-rationality of every $R_\p$ of height two. Proposition~\ref{prop:surface} identifies that condition with $(\mathrm H_2)$.
\end{proof}

Part~(b) has a useful consequence without a characteristic assumption: for a Noetherian Cohen--Macaulay UFD of \emph{any} finite dimension, $R[[t]]$ is a UFD if and only if $(\mathrm H_2)$ holds. This follows from Proposition~\ref{prop:height-two}, because a Cohen--Macaulay ring satisfies $(S_3)$.

For a general nonexcellent $\mathbb Q$-UFD of dimension four, conditions (a) and (b) leave one type of prime untreated: a nonmaximal prime $\p$ of height three with $\depth R_\p=2$. Theorem~\ref{thm:local-intro} proves that $R_\p[[t]]$ is not a UFD, but it does not prove that $R[[t]]$ is not a UFD. The coefficient map $R[[t]]\to R_\p[[t]]$ is not ordinarily a localization map, and the descent in Lemma~\ref{lem:descent} requires finite-length self-Ext supported at a \emph{maximal} ideal. The same issue arises at nonmaximal primes of height at least three in greater dimensions. The criterion is therefore complete through dimension three and for excellent $\mathbb Q$-UFDs in every dimension. For arbitrary nonexcellent $\mathbb Q$-UFDs, the first unresolved case occurs in dimension four.

\begin{proposition}[Conditional reduction to depth-two local rings]
\label{prop:conditional-depth-two}
Assume the following localization implication:
\begin{equation}\tag{$\mathrm{LG}$}\label{eq:LG}
 R[[t]]\text{ is a UFD}
 \quad\Longrightarrow\quad
 R_\p[[t]]\text{ is a UFD}
 \qquad\text{whenever }\depth R_\p=2,
\end{equation}
for every Noetherian $\mathbb Q$-UFD $R$ of finite dimension and every prime
$\p\subset R$. Then, for every such $R$, the following are equivalent:
\begin{enumerate}[label=\textup{(\roman*)},leftmargin=2.7em]
\item $R[[t]]$ is a UFD;
\item $R_\p[[t]]$ is a UFD for every $\p$ with $\depth R_\p=2$;
\item $R$ satisfies $(S_3)$ and $(\mathrm H_2)$.
\end{enumerate}
Thus an all-dimensional proof of \eqref{eq:LG} would make the UFD property of
$R[[t]]$ completely determined by the local depth-two cases. It would also
extend the main criterion of this paper, without excellence or completion
hypotheses, to Noetherian $\mathbb Q$-UFDs of every finite dimension.
\end{proposition}
\begin{proof}
The implication \textup{(i)}$\Rightarrow$\textup{(ii)} is
\eqref{eq:LG}. The converse is Lipman's unrestricted Danilov--Samuel
theorem \cite[Introduction, pp.~203--204]{Lipman1979}.

Corollary~\ref{cor:all-dimensions}\textup{(b)} gives
\textup{(iii)}$\Rightarrow$\textup{(i)}. Conversely, assume
\textup{(i)}. Proposition~\ref{prop:height-two} gives $(\mathrm H_2)$.
If $(S_3)$ failed, then normality of the UFD $R$ and hence $(S_2)$ would
give a prime $\p$ with $\htop\p\geq3$ and $\depth R_\p=2$.
By \eqref{eq:LG}, $R_\p[[t]]$ would be a UFD, whereas
Theorem~\ref{thm:local-intro}, applied to the local $\mathbb Q$-UFD
$R_\p$, says that it is not. This contradiction proves $(S_3)$.
\end{proof}

We record a further criterion that does not assume excellence. A Noetherian ring $T$ satisfies $(S_2)$ if $\depth T_\q\geq\min\{2,\htop\q\}$ at every prime $\q$.

\begin{theorem}[Criterion with $(S_2)$ maximal completions]\label{thm:analytically-s2}
Let $R$ be a Noetherian UFD containing $\mathbb Q$, of arbitrary finite dimension. Suppose that $\widehat{R_\m}$ satisfies $(S_2)$ for every maximal ideal $\m$ of $R$. Then $R[[t]]$ is a UFD if and only if $R$ satisfies $(S_3)$ and $(\mathrm H_2)$.
\end{theorem}
\begin{proof}
It remains to prove that the UFD property of $R[[t]]$ implies $(S_3)$. We first show that $H^2_J(R)$ is a finite $R$-module for every ideal $J$ of height at least three. Fix a maximal ideal $\m\supseteq J$ and write $T=\widehat{R_\m}$. Every prime containing $JT$ has height at least three by flat going down. The complete local ring $T$ is a quotient of a regular complete local ring and is catenary. It also satisfies $(S_2)$; hence it is equidimensional \cite[Lemma~51.3.2, Tag~0FIW]{Stacks}. Suppose $\q\notin V(JT)$ has an immediate specialization $\mathfrak r\in V(JT)$. Catenarity and equidimensionality give $\htop\q=\htop\mathfrak r-1\geq2$, so $\depth T_\q\geq2$. The finiteness theorem of SGA~2 now makes $H^2_{JT}(T)$ finite \cite[Expos\'e VIII, Corollaire~2.3]{SGA2}. Flat base change and faithful flatness of $R_\m\to T$ make $H^2_J(R)_\m$ finite over $R_\m$. This holds at every maximal ideal. Normality of $R$ gives $H^0_J(R)=H^1_J(R)=0$; the local-to-global finiteness theorem for local cohomology now makes $H^2_J(R)$ finite \cite[Lemma~51.7.2, Tag~0AW7]{Stacks}.

If $(S_3)$ failed, normality of $R$ would give a prime $\p$ of height at least three with $\depth R_\p=2$. Then $H^2_\p(R)$ is nonzero, since its localization at $\p$ is $H^2_{\p R_\p}(R_\p)\ne0$. It is finite by the preceding paragraph. Proposition~\ref{prop:finite-H2} contradicts the assumption that $R[[t]]$ is a UFD. The converse is Corollary~\ref{cor:all-dimensions}(b).
\end{proof}

The same criterion applies under another hypothesis that can be checked without passing to completions.

\begin{corollary}[Rings locally embedded in regular schemes]\label{cor:regular-embedding}
Let $R$ be a Noetherian UFD containing $\mathbb Q$ and suppose that $R$ is catenary and $\Spec R$ admits an open cover by schemes that embed as closed subschemes of regular Noetherian schemes. Then $R[[t]]$ is a UFD if and only if $R$ satisfies $(S_3)$ and $(\mathrm H_2)$. In particular, this holds when $R$ is a quotient of a regular Noetherian ring.
\end{corollary}
\begin{proof}
Let $J$ have height at least three. If a prime $\q\notin V(J)$ has an immediate specialization $\mathfrak r\in V(J)$, catenarity of the domain $R$ gives $\htop\q=\htop\mathfrak r-1\geq2$. Normality gives $\depth R_\q\geq2$. By the finiteness theorem of SGA~2, $H^2_J(R)$ is finite \cite[Expos\'e VIII, Corollaire~2.3]{SGA2}. The argument in the last paragraph of the proof of Theorem~\ref{thm:analytically-s2} proves necessity of $(S_3)$; Proposition~\ref{prop:height-two} gives the other implications. A quotient of a regular Noetherian ring is catenary and admits the indicated closed immersion.
\end{proof}

In particular, under either set of hypotheses, the conjectured dimension-four criterion is established. A counterexample to that criterion must have a maximal ideal $\m$ such that $\widehat{R_\m}$ fails $(S_2)$, as well as a nonmaximal height-three prime $\p$ with $\depth R_\p=2$; it also cannot satisfy the hypotheses of Corollary~\ref{cor:regular-embedding}.

The following lemma gives an equivalent test for the completion hypothesis. The generic formal fiber means $T\otimes_A\operatorname{Frac}(A)$.

\begin{lemma}[The low-codimension formal fibers]\label{lem:fiber-test}
Let $(A,\mathfrak n)$ be a normal Noetherian local domain and $T=\widehat A$. Then $T$ satisfies $(S_2)$ if and only if $T\otimes_A\operatorname{Frac}(A)$ satisfies $(S_2)$ and $T\otimes_A\kappa(\mathfrak a)$ satisfies $(S_1)$ for every height-one prime $\mathfrak a\subset A$.
\end{lemma}
\begin{proof}
For $\q\in\Spec T$, put $\mathfrak a=\q\cap A$ and $F=(T\otimes_A\kappa(\mathfrak a))_\q$. The flat local map $A_\mathfrak a\to T_\q$ gives
\[
\depth T_\q=\depth A_\mathfrak a+\depth F,
\qquad \htop\q=\htop\mathfrak a+\dim F
\]
\cite[Lemmas~10.163.2 and~10.112.7, Tags~0337 and~00ON]{Stacks}. If $\htop\mathfrak a\geq2$, normality gives $\depth T_\q\geq2$. If $\mathfrak a=(0)$, the stated $(S_2)$ condition on the generic fiber gives the desired depth. If $\htop\mathfrak a=1$, then $A_\mathfrak a$ is a DVR; the $(S_1)$ condition on its fiber gives $\depth T_\q\geq\min\{2,\htop\q\}$. The same equalities give the converse by subtracting $\depth A_\mathfrak a=0$ or $1$ in the first two cases.
\end{proof}

Thus this formal-fiber test is another formulation of the hypothesis of Theorem~\ref{thm:analytically-s2}. That hypothesis is sufficient, but the following example shows it is not necessary.

\begin{example}[A UFD power series ring with a non-$(S_2)$ completion]\label{ex:nons2-completion}
Let $C=\mathbb Q[[x,y]]$ and
\[
B=\mathbb Q+(x,y)^2C
 =\mathbb Q[[x^2,xy,y^2,x^3,x^2y,xy^2,y^3]]\subset C,
\qquad T=B[[z,w]].
\]
The displayed generators show that $B$ is Noetherian and complete. Its maximal ideal is $\mathfrak n=(x,y)^2C$, and $C/B\simeq\mathbb Qx\oplus\mathbb Qy$ has finite length. The exact sequence $0\to B\to C\to C/B\to0$, together with $\depth_{\mathfrak n}C=2$, gives $\depth B=1$; also $\dim B=2$. Therefore $T$ is a complete local domain with $\dim T=4$ and $\depth T=3$. For the prime $\q=\mathfrak nT$, the closed fibre of the flat local map $B\to T_\q$ is $\operatorname{Frac}(\mathbb Q[[z,w]])$. Thus $\htop\q=2$ and $\depth T_\q=1$. In particular, $T$ fails $(S_2)$ although it is a domain. The finite extension $T\subset C[[z,w]]$ is an isomorphism outside $V(\q)$, so $T_\mathfrak s$ is regular for every height-two prime $\mathfrak s\ne\q$.

Heitmann's isolated-singularity completion theorem \cite[Main Theorem]{HeitmannIsolated} gives a four-dimensional local UFD $(R,\m)$ with residue field $\mathbb Q$ and completion $\widehat R\simeq T$, such that $R_\p$ is regular for every $\p\ne\m$. Its application to a complete ring of depth greater than two is explicit in \cite[proof of Theorem~1.1(1)]{KimuraEtAl}. Every nonzero integer is a unit in $R$, so $\mathbb Q\subset R$. Completion preserves depth: $\depth R=3$, and $R$ satisfies $(S_3)$. Every height-two localization of $R$ is regular, so its power series ring is a UFD. Lipman's Danilov--Samuel theorem gives $R[[t]]$ a UFD \cite[(D--S), introduction, pp.~203--204]{Lipman1979}. Since $T$ is a domain, Ratliff's criterion also makes $R$ universally catenary \cite[Proposition~15.111.5, Tag~0AW6]{Stacks}. This example has no height-three prime of depth two.
\end{example}

Any local UFD $V$ with completion $T$ is regular at its height-two primes. For such a prime $\mathfrak a\subset V$, a prime $\mathfrak s\subset T$ minimal over $\mathfrak aT$ has height two and $\depth T_\mathfrak s\geq\depth V_\mathfrak a=2$. Hence $\mathfrak s\ne\q$, and the regularity of $T_\mathfrak s$ descends to $V_\mathfrak a$. Thus all the UFDs constructed from this completion satisfy $(\mathrm H_2)$.

The isomorphism class of the completion does not by itself determine whether the series ring is a UFD. Indeed, $T_0=B[[z]]$ has dimension three and depth two. Heitmann's UFD-completion theorem gives a three-dimensional local UFD $(A,\mathfrak a)$ with $\widehat A\simeq T_0$ \cite[Theorems~1 and~8]{Heitmann}. Set $R'=(A[w])_{(\mathfrak a,w)}$. Then $R'$ is a four-dimensional local UFD, $\widehat{R'}\simeq T_0[[w]]=T$, and the nonmaximal prime $\mathfrak p=\mathfrak aR'$ has height three and $\depth R'_{\mathfrak p}=2$. Theorem~\ref{thm:local-intro} says that $A[[t]]$ is not a UFD. If $R'[[t]]$ were a UFD, Danilov's faithfully flat descent theorem applied to $A\to R'$ would make $A[[t]]$ a UFD \cite[Proposition~2, p.~370]{Danilov}, a contradiction. Thus the two UFDs $R$ and $R'$ have the same domain as completion, but $R[[t]]$ is a UFD and $R'[[t]]$ is not. Both coefficient rings are universally catenary by Ratliff's criterion.

The following proposition locates the formal-fiber defect required by a dimension-four counterexample.

\begin{proposition}[A necessary configuration in dimension four]\label{prop:four-boundary}
Let $R$ be a four-dimensional Noetherian UFD containing $\mathbb Q$. Suppose that $R[[t]]$ is a UFD and that $\p\subset R$ has height three and $\depth R_\p=2$. Then there are a maximal ideal $\m\supset\p$ and primes $\q\subset\mathfrak r$ in $T=\widehat{R_\m}$ such that
\[
\begin{gathered}
\p T\not\subseteq\q,\quad \p T\subseteq\mathfrak r,\quad
\htop(\mathfrak r/\q)=1,\\
\depth T_\q=1,\quad \dim(T/\q)\in\{1,2\},\quad
\htop(\q\cap R_\m)\leq1.
\end{gathered}
\]
If $T$ is equidimensional, then either its generic formal fiber fails $(S_2)$ at $\q$, or its formal fiber over the height-one prime $\q\cap R_\m$ fails $(S_1)$ at $\q$.
\end{proposition}
\begin{proof}
The module $H^2_\p(R)$ is nonzero, since
\[
H^2_\p(R)_\p=H^2_{\p R_\p}(R_\p)\ne0.
\]
It cannot be finite by Proposition~\ref{prop:finite-H2}. Since $H^0_\p(R)=H^1_\p(R)=0$, the local-to-global finiteness theorem \cite[Lemma~51.7.2, Tag~0AW7]{Stacks} supplies a maximal $\m\supset\p$ for which $H^2_\p(R)_\m$ is not finite. By flat base change and faithful flatness, $H^2_{\p T}(T)$ is not finite. The contrapositive of the SGA~2 finiteness criterion produces $\q\notin V(\p T)$ and an immediate specialization $\mathfrak r\in V(\p T)$ with $\depth T_\q<2$ \cite[Expos\'e VIII, Corollaire~2.3]{SGA2}. Set $\q_0=\q\cap R_\m$. The flat local map $(R_\m)_{\q_0}\to T_\q$ gives $\depth T_\q\geq\depth (R_\m)_{\q_0}$. Normality of $R$ implies $\depth (R_\m)_{\q_0}\geq\min\{2,\htop\q_0\}$, so $\htop\q_0\leq1$.

The prime $\p$ is nonmaximal and has height three, so $\m$ has height four and $\dim R_\m/\p R_\m=1$. Completion preserves this dimension: $\dim T/\p T=1$. Consequently $\dim T/\mathfrak r$ is zero or one. Since $T/\q$ is a complete local domain and $\mathfrak r/\q$ has height one, catenarity gives $\dim T/\q=\dim T/\mathfrak r+1\in\{1,2\}$. The maximal-depth obstruction (Theorem~\ref{thm:local-intro}) implies $\depth R_\m\geq3$, hence $\depth T\geq3$. Every associated prime $\mathfrak s$ of $T$ satisfies $\dim T/\mathfrak s\geq3$ \cite[Lemma~10.72.9, Tag~0BK4]{Stacks}. Thus $\q\notin\Ass T$, and $\depth T_\q\geq1$; the opposite inequality was obtained above.

Finally, suppose $T$ is equidimensional. Since $T$ is complete, $\htop\q=\dim T-\dim T/\q\geq2$. If $\q_0=(0)$, the generic formal fiber has local ring $T_\q$, of depth one and dimension at least two; it fails $(S_2)$. If $\htop\q_0=1$, the flat depth and dimension formulas in Lemma~\ref{lem:fiber-test} show that the formal fiber over $\q_0$ has depth zero and local dimension $\htop\q-1\geq1$ at $\q$; it fails $(S_1)$.
\end{proof}

In positive or mixed characteristic, the depth-three necessity in Theorem~\ref{thm:local-intro} is not established by the exponential and conductor arguments used here. Proposition~\ref{prop:height-two} still gives its stated necessary surface condition and sufficient $(S_3)$-plus-surface condition.

\section{Excellent UFDs in arbitrary characteristic}\label{sec:arbitrary-characteristic}

We conclude with a coefficient-ring criterion in positive and mixed characteristic. Its local Witt-vector component comes from Danilov's Proposition~5 and Lemma~7 \cite[pp.~375--378]{Danilov}. The new ingredient is the passage of the UFD property from $R[[t]]$ to $R_\p[[t]]$.

Let $R$ be a Noetherian UFD of dimension at most three. For $\p\in\Spec R$, set
\[
 A_\p=R_\p,\qquad \m_\p=\p A_\p,\qquad
 U_\p=\Spec A_\p\setminus\{\m_\p\}.
\]
If $\kappa(\p)=A_\p/\m_\p$ has prime characteristic $\ell>0$, let $W_N\mathcal O_{U_\p}$ denote the sheaf of length-$N$ $\ell$-typical Witt vectors on $U_\p$. Write
\begin{equation}\label{eq:witt-image}
 Z_N(A_\p)=\operatorname{im}\!\left[
 H^1(U_\p,W_N\mathcal O_{U_\p})
 \xrightarrow{\text{first Witt coordinate}}
 H^1(U_\p,\mathcal O_{U_\p})\right].
\end{equation}
This definition also applies when $\operatorname{char}A_\p=0$ and $\operatorname{char}\kappa(\p)=\ell$.

\begin{theorem}[Criterion under excellence]\label{thm:arbitrary-characteristic}
Let $R$ be an excellent Noetherian UFD with $\dim R\leq3$. Then $R[[t]]$ is a UFD if and only if both of the following conditions hold:
\begin{enumerate}[label=\textup{(\roman*)},leftmargin=2.7em]
\item for every height-two prime $\p$ and all $u,v\in A_\p$ satisfying $\sqrt{(u,v)A_\p}=\m_\p$, one has $\overline{((u,v)A_\p)^2}\subseteq(u,v)A_\p$;
\item for every height-three prime $\p$ with $\depth A_\p=2$, the field $\kappa(\p)$ has characteristic $\ell>0$ and $Z_N(A_\p)=0$ for some integer $N\geq1$.
\end{enumerate}
The integer $N$ may depend on $\p$. No condition is imposed at a height-three prime where $\depth A_\p=3$.
\end{theorem}

Condition~\textup{(ii)} can be checked in the coefficient ring. If $\m_\p=(f_1,\ldots,f_r)$, the equality $Z_N(A_\p)=0$ says that, for every family $c_{ij}\in W_N((A_\p)_{f_if_j})$ satisfying $c_{ij}+c_{jk}=c_{ik}$ under Witt addition, there are $b_i\in(A_\p)_{f_i}$ such that the first coordinate of $c_{ij}$ equals $b_j-b_i$ for all $i,j$. Thus no power-series ring appears among the hypotheses of Theorem~\ref{thm:arbitrary-characteristic}.

We first prove the localization step. Its excellence hypothesis concerns only the local ring $R_\p$.

\begin{lemma}[Finite-jet detection]\label{lem:char-jet}
Let $(A,\m)$ be an excellent normal Noetherian local domain, let $T=A[[t]]$, and let $U=\operatorname{Reg}(\Spec A)$. Put $W=U\times_{\Spec A}\Spec T$ and $U_n=U\times_{\Spec A}\Spec(A[t]/t^n)$ for $n\geq1$. Then restriction identifies $\Cl(T)$ with $\Pic(W)$, and the map
\[
 \Pic(W)\longrightarrow\varprojlim_{n\geq1}\Pic(U_n)
\]
is injective.
\end{lemma}
\begin{proof}
The map $A[t]\to T$ is the $(t)$-adic completion of the excellent ring $A[t]$, hence is regular \cite[Lemma~15.51.14, Tag~0AH2]{Stacks}. The regular map $A\to T$ makes $W$ regular. Since $A$ is normal and $A\to T$ is flat, the complement of $W$ in $\Spec T$ contains no height-one prime. Hence $\Cl(T)=\Cl(W)=\Pic(W)$.

As $U$ contains every height-one point of $\Spec A$, normality gives $\Gamma(U,\mathcal O_U)=A$, and therefore $\Gamma(U_n,\mathcal O_{U_n})=A[t]/t^n$. The reduction maps between the unit groups of these rings are surjective. Suppose $L\in\Pic(W)$ restricts to the trivial bundle on every $U_n$. Choose compatible trivializations, and let $M=\Gamma(W,L)$, a finite rank-one reflexive $T$-module. The trivializations give a homomorphism
\[
 \varphi:M\longrightarrow\varprojlim_n\Gamma(U_n,\mathcal O_{U_n})=T.
\]
At a prime $P\in W$ containing $t$, its reduction modulo $t$ is the chosen trivialization on $U_1$; Nakayama's lemma makes $\varphi_P$ an isomorphism. If $Q\subset T$ has height one and does not contain $t$, take $P$ minimal over $Q+(t)$. The ring $T$ is local, so this ideal is proper. The principal ideal theorem gives $\htop P=2$, and $P/(t)$ has height one in $A$. Thus $P\in W$ and $\varphi_Q$ is an isomorphism as well. Two reflexive $T$-modules with a map that is an isomorphism at every height-one prime are isomorphic. Consequently $M\simeq T$ and $L$ is trivial.
\end{proof}

\begin{lemma}[Rescaling a pointed \`etale neighborhood]\label{lem:char-rescaling}
Let $R$ be a domain, $\p\in\Spec R$, $s\in R\setminus\p$, and let $E$ be an affine \`etale $R_s[t]$-algebra of finite presentation with a specified isomorphism $E/tE\simeq R_s$. Set $A=R_\p$, and let $\phi:E\otimes_{R_s}A\to A[[t]]$ be the unique $A[t]$-algebra map reducing to the specified map at $t=0$. After replacing $s$ by a multiple $s'\in R\setminus\p$, there are $m\geq1$ and a map
\[
 \theta:E\otimes_{R_s}R_{s'}\longrightarrow(R[[t]])_{s'},
 \qquad \theta(t)=(s')^m t,
\]
whose composition with $(R[[t]])_{s'}\to A[[t]]$ equals $\phi$ followed by the automorphism $t\mapsto(s')^m t$ of $A[[t]]$.
\end{lemma}
\begin{proof}
The special fiber of $E\otimes_{R_s}A$ is $\Spec A$. Its closed point is a specialization of every point on that fiber. Thus an open chart about that point contains the entire special fiber. A standard \`etale chart is available by \cite[Proposition~10.144.4, Tag~00UE]{Stacks}. Spread the chart to $R_{s'}$ and replace $s$ by $s'$; we may write it as
\[
\begin{gathered}
E=R_s[t,z,1/g(t,z)]/(f(t,z)),\qquad f(0,a)=0,\\
j=f_z(0,a)\in R_s^\times,\qquad g_0=g(0,a)\in R_s^\times.
\end{gathered}
\]
for some $a\in R_s$. Including $f_z$ among the inverted elements, if necessary, makes the displayed chart \`etale. Since its special fiber contains the closed point of $\Spec A$, we may enlarge $s$ so that $j$ and $g_0$ are units of $R_s$.

Normalize the equation as
\[
 j^{-1}f(t,a+Y)=Y+
 \sum_{\substack{i,k\geq0\\i\geq1\ \mathrm{or}\ k\geq2}}
 c_{ik}t^iY^k,\qquad c_{ik}\in R_s.
\]
Choose $d\geq1$ so that $s^dc_{ik}\in R$ for all the finitely many coefficients, and the corresponding normalized coefficients of $g(t,a+Y)-g_0$ have the same property. Then
\[
 P(t,Z)=s^{-d}j^{-1}f(s^{2d}t,a+s^dZ)
 =Z+\sum_{i,k}c_{ik}s^{(2i+k-1)d}t^iZ^k
 \in R[t,Z].
\]
It satisfies $P(0,0)=0$ and $P_Z(0,0)=1$, so recursion determines a unique $Z(t)\in tR[[t]]$ with $P(t,Z(t))=0$. Send $z$ to $a+s^dZ(t)$ and the base variable $t$ to $s^{2d}t$. After dividing by $g_0$, the image of $g$ belongs to $1+tR[[t]]$ and is invertible in $(R[[t]])_s$. This constructs $\theta$, with $m=2d$. Its asserted behavior over $A$ follows from the uniqueness of the lift across every square-zero quotient $A[t]/t^{n+1}\to A[t]/t^n$ for an \`etale algebra.
\end{proof}

\begin{proposition}[Coefficient localization under excellence]\label{prop:char-localization}
Let $R$ be a Noetherian UFD and $\p\in\Spec R$. If $R_\p$ is excellent and $R[[t]]$ is a UFD, then $R_\p[[t]]$ is a UFD.
\end{proposition}
\begin{proof}
Put $A=R_\p$, $T=A[[t]]$, and $U=\operatorname{Reg}(\Spec A)$. Suppose $T$ is not a UFD. Lemma~\ref{lem:char-jet} gives $n\geq1$ and a line bundle $L$ on $U\times_A\Spec T$ whose restriction to $U_n$ is nontrivial.

Let $H$ be the henselization of the pair $(A[t],(t))$. It is Noetherian, its $(t)$-adic completion is $T$, and $H/t^nH=A[t]/t^n$ \cite[Lemma~15.12.4, Tag~0AGV]{Stacks}. A finite principal-open cover of $U$ expresses the transition functions of $L$, their inverses, and their cocycle identities by finitely many polynomial equations over $H$, after denominators are cleared. Approximation for the henselian G-ring pair $(H,(t))$ produces a line bundle on $U\times_A\Spec H$ whose restriction to $U_n$ is isomorphic to $L|_{U_n}$ \cite[Lemma~16.14.1, Tag~0AH5]{Stacks}. The ring $H$ is a filtered colimit of affine pointed \`etale $A[t]$-algebras $E$ with $E/tE\simeq A$. Since the cover and transition data are finite, our bundle descends to one such $E$. Spreading out from $A=R_\p$, after passing to $R_s$ with $s\notin\p$, gives an \`etale $R_s[t]$-algebra $E_s$ with $E_s/tE_s\simeq R_s$, a quasi-compact open $V_s\subset\Spec R_s$ whose base change to $A$ is $U$, and a line bundle on $V_s\times_{\Spec R_s}\Spec E_s$.

Apply Lemma~\ref{lem:char-rescaling}, enlarging $s$ as needed, and pull this bundle back to $V_s\times_{\Spec R_s}\Spec (R[[t]])_s$. This is an open subscheme of the spectrum of a UFD, so its Picard group is zero. Its further pullback to $U_n$ is $L|_{U_n}$ after the automorphism $t\mapsto s^m t$ of $A[t]/t^n$; it is nontrivial. This contradiction proves that $T$ is a UFD.
\end{proof}

\begin{proof}[Proof of Theorem~\ref{thm:arbitrary-characteristic}]
Suppose $R[[t]]$ is a UFD. Proposition~\ref{prop:height-two} proves condition~\textup{(i)} without a characteristic assumption. By Proposition~\ref{prop:char-localization}, $A_\p[[t]]$ is a UFD for every height-three prime $\p$. Fix such a prime with $\depth A_\p=2$. Danilov's Proposition~5(1) makes its infinite Witt obstruction vanish \cite[pp.~375--376]{Danilov}. If $\operatorname{char}\kappa(\p)=0$, the obstruction is $H^1(U_\p,\mathcal O_{U_\p})=H^2_{\m_\p}(A_\p)\ne0$, a contradiction. Thus $\operatorname{char}\kappa(\p)=\ell>0$. Danilov's Lemma~7 identifies the vanishing of the obstruction with $\bigcap_{N\geq1}Z_N(A_\p)=0$ \cite[p.~377]{Danilov}. The $Z_N(A_\p)$ form a descending chain of $A_\p$-submodules of the Artinian module $H^2_{\m_\p}(A_\p)$. The chain stabilizes, so $Z_N(A_\p)=0$ for some $N$. This proves condition~\textup{(ii)}.

Conversely, assume~\textup{(i)} and~\textup{(ii)}. For every prime of height two, Proposition~\ref{prop:surface} makes $A_\p[[t]]$ a UFD. The same is true at primes of height at most one because the coefficient local rings are regular. Let $\p$ have height three. If $\depth A_\p=2$, all its nonmaximal localizations have the discrete class-group property by the preceding cases; condition~\textup{(ii)}, Danilov's Lemma~7 and Proposition~5(2) then make $A_\p[[t]]$ a UFD \cite[pp.~375--378]{Danilov}. If $\depth A_\p=3$, the Danilov--Samuel reduction applies because every localization of $A_\p$ of depth two has height two \cite[(D--S), introduction, pp.~203--204]{Lipman1979}. Thus every maximal localization of $R$ has the discrete class-group property, and Danilov's localness theorem gives $R[[t]]$ a UFD \cite[Proposition~3, p.~370]{Danilov}.
\end{proof}

\begin{remark}[Scope without excellence]
For a three-dimensional local Noetherian UFD $A$ whose completion is normal, Danilov's Proposition~5 also applies without excellence. The height-two test in Theorem~\ref{thm:arbitrary-characteristic}\textup{(i)}, together with the condition at the maximal ideal in~\textup{(ii)}, is then necessary and sufficient for $A[[t]]$ to be a UFD. In dimension at most two, Theorem~\ref{thm:arbitrary-characteristic}\textup{(i)} alone is necessary and sufficient for every Noetherian UFD, without excellence, by Proposition~\ref{prop:height-two}. For arbitrary nonexcellent three-dimensional UFDs, this section does not assert the global equivalence: Proposition~\ref{prop:char-localization} requires excellence of the local ring $R_\p$.
\end{remark}

\end{document}